\documentclass{amsart}

\usepackage{amsfonts}
\usepackage{graphicx}
\usepackage{epstopdf}

\newcommand{\dpl}{\displaystyle}

\newcommand{\beq}{\begin{equation}}
\newcommand{\eeq}{\end{equation}}

\begin{document}

\title[Yang-Baxter maps and Poisson Lie groups]{Poisson Lie groups and Hamiltonian theory of the Yang-Baxter maps}
       \author{Nicolai Reshetikhin}
       \address{Department of  Mathematics, University of California, Berkeley,
CA 94720, USA}
           \email{reshetik@math.berkeley.edu}

       \author{Alexander Veselov}
       \address{Department of Mathematical Sciences,
        Loughborough University, Loughborough,
        Leicestershire, LE11 3TU, UK and
       Landau Institute for Theoretical Physics, Moscow, Russia}
       \email{A.P.Veselov@lboro.ac.uk}

\maketitle




{\small  {\bf Abstract.} We show how  the theory of Poisson Lie groups can be used to establish the Poisson properties of the Yang-Baxter maps and related transfer dynamics. As an example we present the Hamiltonian structure for the matrix KdV soliton interaction.}

\bigskip

\section*{Introduction}
The general theory of the Yang-Baxter maps (also known as the set-theoretical
solutions to the quantum Yang-Baxter equation) was initiated by Drinfeld's question \cite{D}, but the first interesting examples appeared already in the paper \cite{Skl} by Sklyanin.  Important construction of Yang-Baxter maps was derived by Weinstein and Xu \cite{WX}
from the theory of Poisson Lie groups.

The dynamical aspects of this problem in relation with the theory of integrable systems were recently investigated in \cite{V, GV, SV} (see also earlier works \cite{R1,BBR}). 
The general idea is that the corresponding transfer-maps should be integrable in some sense. For the symplectic maps one can use the standard notion of Liouville integrability, so the main question is whether the transfer-maps are symplectic and what is the corresponding Poisson structure.

In this paper we show that the theory of Poisson Lie groups can give the answer to this question
and show how it works for the concrete examples of Yang-Baxter maps from \cite{V, GV}.

\section{Yang-Baxter maps and related transfer-dynamics.}

\subsection{} Let $X$ be a set and $R$ be a map: $$R: X \times X \rightarrow X
\times X.$$ Let $R_{ij}: X^{n} \rightarrow X^{n}, \quad X^{n} = X
\times X \times .....\times X$ be the maps which acts as $R$ on
$i$-th and $j$-th factors and identically on the others. If $P:
X^2 \rightarrow X^2$ is the permutation: $P(x,y) = (y,x)$, then
$$R_{21} = P R P.$$
 $R$ is called the {\it Yang-Baxter map}
if it satisfies the Yang-Baxter relation
\begin{equation}
\label{YB} R_{12} R_{13} R_{23} = R_{23} R_{13} R_{12},
\end{equation}
considered as the equality of the maps of $X \times X \times X$
into itself.  We will call $R$ {\it reversible Yang-Baxter map}.
If additionally $R$ satisfies the relation
\begin{equation}
\label{U}
R_{21} R = Id,
\end{equation}

We will consider also the parameter-dependent Yang-Baxter
maps $R(\lambda,\mu), \lambda, \mu \in {\bf C}$ satisfying the
corresponding version of Yang-Baxter relation
\begin{equation}
\label{sYB} R_{12}(\lambda_1, \lambda_2) R_{13}(\lambda_1,
\lambda_3) R_{23} (\lambda_2, \lambda_3) = R_{23}(\lambda_2,
\lambda_3) R_{13}(\lambda_1, \lambda_3) R_{12}(\lambda_1,
\lambda_2).
\end{equation}
Although this case can be considered as a particular case of the
previous one by introducing $\tilde X = X \times {\bf C}$ and
$\tilde R (x,\lambda; y, \mu) = R (\lambda,\mu) (x,y)$ it is often 
convenient to keep the parameter separately.

Let us define the {\it transfer-maps} $T_i^{(N)}, i = 1,\dots, N$ as the maps of $X^{N}$ into 
itself given by the following formulas:
\begin{equation}
\label{T}
T_i^{(N)} = R_{i i+N-1} R_{i i+N-2} \dots R_{i i+1},
\end{equation}
where the indices are considered modulo $n.$
In particular $T_1^{(N)} = R_{1 N} R_{1 N-1} \dots R_{12}.$ 
 
It is easy to show that for any reversible Yang-Baxter map $R$ the transfer-maps $T_i^{(N)},\quad i =1,\dots, N$ commute with each other:
\begin{equation}
\label{comm}
T_i^{(N)} T_j^{(N)} = T_j^{(N)} T_i^{(N)}
\end{equation}
and satisfy the property
\begin{equation}
\label{prod}
T_1^{(N)} T_2^{(N)} \dots T_N^{(N)} = Id
\end{equation}
and moreover these properties can be used as characteristic for the reversible YB maps
(see \cite{V}). Thus for any reversible YB maps one can define the action of ${\bf Z}^{N-1}$
on $X^{N}$ generated by $T_i^{(N)}$, which we will call {\it transfer-dynamics}.

\subsection{Lax matrices.} Let $R(\lambda,\mu)$ be a parameter-dependent Yang-Baxter
map,  then by its {\it Lax matrix} (or Lax representation) we will mean the matrix ${\mathcal L}(x,\lambda; \zeta)$ depending on the point
$x \in X$, parameter $\lambda$ and additional ("spectral")
parameter $\zeta \in {\bf C},$  which satisfies the following
relation
\begin{equation}
\label{Lax}
 {\mathcal L}(x,\lambda; \zeta){\mathcal L}(y,\mu; \zeta) =  {\mathcal L}(\tilde{y},\mu; \zeta){\mathcal L}(\tilde{x},\lambda; \zeta)
\end{equation}
where $(\tilde{x},\tilde{y}) = R (\lambda,\mu) (x,y).$ It is easy to show
that such a matrix allows to produce the integrals
for the dynamics of the related transfer-maps: the spectrum of the  {\it monodromy matrix} 
\begin{equation}
M =  {\mathcal L}(x_1,\lambda_1, \zeta) {\mathcal L}(x_{2},\lambda_{2},\zeta) \ldots {\mathcal L}(x_n,\lambda_n, \zeta)
\end{equation}
is preserved under the corresponding transfer-maps $T_i^{(n)}, \quad i=1,\ldots, n.$
In \cite{SV} it was shown that for certain class of Yang-Baxter maps one can write down the Lax matrices straight from the map itself.
Namely, suppose that on the set $X$ we have an effective action of a linear group $G
=  GL_N$ and the Yang-Baxter map $R(\lambda, \mu)$ has the
following special form:
\begin{equation}
\label{map} \tilde{x} = {\mathcal K}(y,\mu, \lambda)[x], \quad  \tilde{y} =
{\mathcal L}(x,\lambda, \mu)[y],
\end{equation}
where ${\mathcal K},{\mathcal L} : X \times {\bf C} \times {\bf C} \rightarrow GL_N$ are
some matrix valued functions on $X$ depending on parameters
$\lambda$ and $\mu$ and $A[x]$ denotes the action of $A \in G$ on
$x \in X.$ Then the claim is that ${\mathcal L}(x,\lambda, \zeta)$ is a Lax matrix for
$R.$ The proof easily follows from the Yang-Baxter relation.
Indeed apply both sides of the relation (\ref{sYB}) to a triple $(x_1,x_2,x_3)$ of points in the set $X$
and compare the results for $x_3.$ Then the left-hand side gives us 
$ {\mathcal L}(x_1,\lambda_1; \lambda_3){\mathcal L}(x_2,\lambda_2; \lambda_3)[x_3]$
while the right-hand side is $ {\mathcal L}(\tilde{x}_2,\lambda_2; \lambda_3){\mathcal L}(\tilde{x}_1,\lambda_1; \lambda_3)[x_3],$ which implies the Lax relation (\ref{Lax}) with $\zeta = \lambda_3.$.

Our main observation here is that all examples of the Lax matrices considered in \cite{SV} 
can be naturally interpreted in terms of symplectic geometry of the Poisson Lie groups.

\section{Poisson Lie groups and Yang-Baxter correspondences.}

The following general construction of the Yang-Baxter correspondences
from the Poisson Lie groups can be traced back to the paper \cite{WX} by Weinstein and Xu.

Recall that a Poisson structure on a manifold is a bi-vector field
defining a Lie bracket on the space of smooth functions.

The {\it correspondence} (relation) $\Phi$ from a set $X$ into
itself is a multivalued, partially defined map determined by its
graph $\Gamma_{\Phi}$, which is a subset of $X \times X.$ Now let
$X = M$ be a smooth manifold with Poisson structure $J$ and the
graph $\Gamma_{\Phi}$ is a submanifold of $M \times M.$ Let us
supply the manifold $M \times M$ with the Poisson structure $J
\oplus (- J).$ Recall now that a submanifold $N$ of a Poisson
manifold $P$ is called { \it coisotropic} if the space of
functions on $P$ which vanish on $N$ is closed under Poisson
bracket. The correspondence $\Phi$ is called {\it Poisson} if its
graph $\Gamma_{\Phi}$ is coisotropic submanifold of $M \times M.$
For the usual maps this equivalent to the standard notion of
Poisson map. If the appropriate smoothness condition holds then
the composition of two Poisson correspondences is again Poisson.

For details about Poisson correspondences see \cite{W}.

\subsection{Yang-Baxter correspondences}
Note that for the correspondence $R: X \times X \rightarrow X \times X$ the map $R_{12}$
is defined by its graph which is a subset of $X^{(6)}$ of the form
$\Gamma_R \times \Delta,$
where $\Delta = (x,x), x \in X$ is the diagonal in $X \times X.$ Similarly one can define
$R_{ij}$ and thus {\it Yang-Baxter correspondence} as such $R$ which satisfies the relation
(\ref{YB}). Reversibility is a bit more tricky. The problem is with the notion of the
inverse correspondence $\Phi^{-1}$: the standard formula $\Phi \Phi^{-1} = Id$ seems
to be too restrictive, so we should avoid it.
There are two involutions on $X \times X \times X \times X$:
$$ \pi:  (x_1, x_2; y_1, y_2) = (x_2, x_1; y_2, y_1)$$
and
$$\tau:  (x_1, x_2; y_1, y_2) =  (y_1, y_2; x_1, x_2).$$
We will call the correspondence $R$ {\it reversible}
if $\pi \Gamma_{R} = \tau \Gamma_{R},$ or equivalently
its graph $\Gamma_{R}$ is invariant under the involution
$$\sigma = \pi \tau: (x_1, x_2; y_1, y_2) = (y_2, y_1; x_2, x_1).$$
For the maps this is equivalent to the definition given above.

\subsection{Poisson Lie groups}
A Poisson Lie group is a Lie group with the Poisson structure such
that the group multiplication is a Poisson mapping (see \cite{D2}).

For any Poisson Lie group $G$ one can construct the following
correspondence $\Phi_G: G \times G \rightarrow G \times G$: its
graph $\Gamma_{\Phi}$ is the set of $(g_1, g_2; h_1, h_2) \in G
\times G \times G \times G$ such that
\begin{equation}
g_1 g_2 = h_2 h_1.
\label{fact}
\end{equation}
From the definition of the Poisson Lie group it follows that this correspondence is Poisson (cf. \cite{WX}).

{\bf Proposition 1.}  {\it For any Poisson Lie group $G$ the Poisson correspondence $\Phi_G$ is a reversible Yang-Baxter correspondence.}

This follows from the associativity of multiplication in $G$ (cf. \cite{V}): both left and right sides are the correspondences with the same graph given by the relation $$g_1 g_2 g_3 = h_3 h_2 h_1.$$ Reversibility is obvious.

Since we are primarily interested in genuine maps but not correspondences this construction is too general.  A natural idea therefore is to restrict it to the symplectic leaves in $G$ in order to make the re-factorisation relation
(\ref{fact}) uniquely solvable.

Suppose that there exist a one-parameter family of embeddings of the set $X$ as a symplectic leaf in the Poisson Lie group $G$: $\phi_{\lambda} : X \rightarrow G$ and define the correspondence 
$R(\lambda, \mu): X \times X \rightarrow X \times X$ by the relation
\begin{equation}
\phi_{\lambda}(x) \phi_{\mu}(y) =  \phi_{\mu}(\tilde{y}) \phi_{\lambda}(\tilde{x}).
\label{factl}
\end{equation}

Note that the embedding $\phi_{\lambda} $ induce on $X$ the symplectic structure
\begin{equation}
\label{symp}
\omega_{\lambda} = \phi_{\lambda}^ *(\omega_{G}),
\end{equation}
where $\omega_{G}$ is the symplectic structure on the corresponding leaf in $G$.
Let us introduce on $X \times X$ a symplectic structure as the direct sum $\omega_{\lambda} \oplus \omega_{\mu}.$ In a similar way we define
the symplectic structure $\Omega^{(N)}$ on $X^{(N)}$ as the direct sum of the forms $\omega_i = \omega_{\lambda_i}:$
\begin{equation}
\label{Omega}
\Omega^{(N)} = \omega_{1} \oplus \omega_{2} \oplus ...\oplus \omega_{N}.
\end{equation}

{\bf Theorem 1.} {\it  The correspondence $R(\lambda, \mu)$ defined by (\ref{factl})  is a reversible Yang-Baxter Poisson correspondence.
The related transfer-maps (\ref{T}) are Poisson correspondences with respect to the symplectic structure}  $\Omega^{(N)}.$

This construction is already good enough to produce many interesting examples of symplectic Yang-Baxter maps, including the interaction of matrix solitons \cite{V,GV}. 
To show this we consider the standard Poisson Lie loop group $G = LGL(n, \bf{R})$
and the symplectic leaves of a special form.

\section{Symplectic leaves in $LGL(n, \bf{R})$ and examples of Yang-Baxter maps}

The standard Poisson structure on the loop group $G = LGL(n, \bf{R})$ is given by the Sklyanin formula
\begin{equation}
\{L_1(\zeta),\otimes L_2(\eta) \} = [ \frac{P_{12}}{\zeta - \eta}, L_1(\zeta) \otimes L_2(\eta)],
\label{Skl}
\end{equation}
where $P_{12}$ is the permutation matrix: $P(x \otimes y) = y \otimes x.$

We will consider the symplectic leaves of the form $L(\zeta) = A + B \zeta.$
The relation (\ref{Skl}) is equivalent to the following set of relations
\begin{equation}
\label{leaves}
\{A_1, \otimes A_2\} = P_{12}(A_1 B_2 - B_1 A_2),  \quad \{A_1, \otimes B_2\} = 0, \quad \{B_1, \otimes B_2\} = 0.
\end{equation}

From these relations it follows that the matrix elements of $B$ are in the centre of the corresponding Poisson algebra, so $B$ can be fixed. 

\subsection{Non-degenerate $B$: matrix soliton interaction}
If $B$ is non-degenerate one can assume that $B = I$ is identity matrix. The rest of the relations (\ref{leaves}) determine the standard Lie-Poisson structure on the matrix group $GL(n, \bf{R})$ considered as an open set in the $n \times n$-matrix algebra $M_n$ which can be identified as a space with 
$gl^*(n, \bf{R})$.

Let ${\mathcal O}_k$ be the special symplectic leaf in $gl^*(n, \bf{R})$ consisting of involutions $S = 2P - I,$ where $P$ is a projector of rank $k$: $P^2 = P, S^2 = I.$ Such a projector $P$ is determined by its kernel $K = Ker P$ and its image $L = Im P,$  which are two subspaces of dimension $k$ and $n-k$ respectively such that $K \oplus L = \bf{R}^n.$ Thus $X$ is an open subset of the product of two Grassmannians $G(n, k)$ and $G(n, n-k)$. Note that the Grassmannian $G(n, n-k)$ is naturally isomorphic to the Grassmannian $G(n, k)$ in the dual space.

Let $X = {\mathcal O}_k$ and define the embedding $\phi_{\lambda}$ of $X$ into $G$ by the formula
\begin{equation}
\label{phi}
\phi_{\lambda}(S) = \zeta I + \lambda S = (\zeta - \lambda) I + 2 \lambda P .
\end{equation}

The claim is that in this case our construction leads to the Yang-Baxter maps
introduced in \cite{GV} in relation with matrix KdV equation. For this it is enough to compare
the formula (13) of the paper \cite{GV} with our formulas (\ref{factl}), (\ref{phi}).

In terms of the subspaces $K$ and $L$ the corresponding Yang-Baxter map $R(\lambda_1, \lambda_2)$ is given by the following formulas:
\begin{equation}
\label{K1}
\tilde K _1 = (I - \frac{2 \lambda_2}{\lambda_1
+\lambda_2}P_2) K_1,
\end{equation}
\begin{equation}
\label{K2}\tilde K_2 = (I - \frac{2 \lambda_1}{\lambda_1
+\lambda_2}P_1) K_2,
\end{equation}
\begin{equation}
\label{L1} \tilde L_1 = (I + \frac{2 \lambda_2}{\lambda_1
-\lambda_2}P_2) L_1,
\end{equation}
\begin{equation}
\label{L2} \tilde L_2 = (I + \frac{2 \lambda_1}{\lambda_2
-\lambda_1}P_1) L_2,
\end{equation}
which give the solution of the corresponding re-factorisation problem (\ref{factl}),
which is uniquely solvable in this case (see \cite{GV}). 

Let us discuss the rank $k=1$ case in more detail.

Let $V$ be an $n$-dimensional real vector space, $V^*$ be its dual,
$<p,q>$ denote the canonical pairing between vector $q \in V$ and covector $p \in V^*.$
The projectors in $V$ have the form $P = p \otimes q,$ where covector $p$ and vector $q$
satisfy the relation $<p,q>=1.$ Obviously such representation is not unique: we can replace
$p \rightarrow e^t p, q \rightarrow e^{-t} q.$ Since this is the Hamiltonian action with the Hamiltonian $H = <p,q>$  the quotient symplectic space ${\mathcal O}_1$ is the result of the Hamiltonian reduction of $T^* V$ with the standard structure $dp \wedge dq$ on the level $H = 1.$ One can check that the symplectic structure $\omega_{\lambda}$ on the space of projectors induced by the embedding (\ref{phi}) can be described as the result of the same Hamiltonian reduction on the level $H = \lambda^2.$

The map $R(\lambda_1, \lambda_2): (p_1, q_1; p_2, q_2)
\rightarrow (\tilde{p_1}, \tilde{q_1}; \tilde{p_2},
\tilde{q_2)}$ in this case has the form 
\beq \label{fi} \tilde{p_1} = p_1+\frac{\dpl
2\lambda_2(p_1,q_2)}{\dpl (\lambda_1-
\lambda_2)(p_2,q_2)}p_2, \qquad \tilde{q_1} =
q_{1}+\frac{\dpl 2\lambda_2(p_2,q_1)}{\dpl (\lambda_1-
\lambda_2)(p_2,q_2)}q_2, \eeq

\beq \label{fj} \tilde{p_2} = p_2+\frac{\dpl
2\lambda_1(p_2,q_1)}{\dpl (\lambda_2-
\lambda_1)(p_1,q_1)}p_1,\qquad \tilde{q_2} =
q_2+\frac{\dpl 2\lambda_1(p_1,q_2)}{\dpl (\lambda_2
-\lambda_1)(p_1,q_1)}q_1\eeq
and describes the interaction of two matrix solitons
for the matrix KdV equation 
$$U_t+3UU_x+3U_xU+U_{xxx}=0 $$
(see for the details \cite{GV}).

{\bf Theorem 2.}  {\it The transfer-dynamics related to the interaction of matrix solitons is symplectic with respect to the corresponding form (\ref{Omega}). }

The coefficients of the characteristic polynomial $\chi(\zeta, \eta) = \det ({\mathcal L} - \eta I)$
of the monodromy matrix 
\begin{equation}
\label{mon}
M(\zeta) =  {\mathcal L}(P_1,\lambda_1, \zeta) {\mathcal L}(P_{2},\lambda_{2},\zeta) \ldots {\mathcal L}(P_N,\lambda_N, \zeta),
\end{equation}
where ${\mathcal L}(P, \lambda, \zeta) = (\zeta - \lambda) I + 2 \lambda P,$ give a set of involutive integrals for transfer dynamics. However already the case $N=2$ shows that this is not enough for complete integrability, so we need more integrals. For example, it is easy to show that all the entries of the matrix $J = \sum_{i=1}^N \lambda_i P_i$ are the integrals as well.
This corresponds to the matrix-valued integral $I _0= \int U dx$ for the matrix KdV equation and is related to the symmetry of the system under the diagonal coadjoint action of $GL(n)$.
We leave the complete investigation of the corresponding dynamics for the future only mentioning here the paper \cite{R}, where a related question was investigated in full generality for the simple Lie groups with standard Poisson Lie structure.

\bigskip

Let us discuss now the symplectic leaves with degenerate $B.$

\subsection{Degenerate $B$:  Adler's maps as Hamiltonian reduction}

Probably the simplest example of the Yang-Baxter maps has the form
\begin{equation}
\label{Adler} \tilde{x} = y - \frac{\lambda - \mu}{x +y} \qquad
\tilde{y} = x - \frac{\mu - \lambda}{x+y},
\end{equation}
where $x, y \in X = {\bf R}.$
This map (modulo additional permutation) first appeared in Adler's
paper \cite{A} as a symmetry of the periodic dressing chain \cite
{VS}.  It was actually derived from the re-factorization problem for the matrix $$ {\mathcal L} =
\Bigl(
\begin{array}{cc}
x &  x^2 + \lambda - \zeta\\ 1 & x
\end{array}
\Bigl), $$
so the Lax pair and the integrals for this map was known from the very beginning.

What we would like to explain here is the Hamiltonian properties of the corresponding transfer dynamics
using the Poisson Lie groups. We are going to show that our construction combined with Hamiltonian reduction naturally leads to the Poisson structure of the dressing chain.

Recall that {\it periodic dressing chain} is the following system of ordinary differential equations 
\begin{equation}
\label{dc}
(f_i + f_{i+1})' = f_i^2  - f_{i+1}^2 + \lambda_i - \lambda_{i +1},\,\, i=1,2,...,N,
\end{equation}
$\lambda_1, \lambda_2,..., \lambda_N$ are some constant parameters 
and we assume that $f_{N+1} = f_1, \lambda_{N+1} = \lambda_1$.
This system can be resolved with respect to $f'_i$ only if the number $N$ is odd.
In that case this is a Hamiltonian system with the Hamiltonian 
$$H = \sum_{i =1}^{N} ( \frac{1}{3} f_i^3 + \lambda_i f_i)$$
with respect to the following Poisson structure \cite{VS}.
Let us introduce the variables $g_i = f_i + f_{i+1}, i = 1,..., N,$
then $$\{g_i, g_{i+1}\} = 1 = - \{g_{i+1}, g_i\} , \,\,i = 1, ..., N$$
and all other Poisson brackets are zero.
For the original variables $f_i$ we have
$$\{f_i, f_j\} = (-1)^{j-i+1 (mod N)}$$
if $i \neq j$ and $0$ otherwise.

To explain the relation of this Poisson structure with Poisson Lie groups
let us consider first the symplectic leaves in $LGL(2, {\bf R})$ of the form
$L = A + B \zeta$ with 
$$ B = \Bigl(
\begin{array}{cc}
0 &  1 \\ 0  & 0
\end{array}
\Bigl). $$ 
Writing $A$ as
$$ A = \Bigl(
\begin{array}{cc}
a &  b \\ c  & d
\end{array}
\Bigl) $$ 
we have the following Poisson brackets for the matrix elements from the relation
(\ref{Skl}):
$$
\label{abcd} 
\{a, b\} = a, \quad \{a, d\} = c, \quad \{a, c\} = 0, \quad \{b, c\} = 0, \quad \{b, d\} = d, \quad
\{c, d\} = 0. 
$$
It is easy to check that the centre of this Poisson algebra is generated by two functions: $C_1 = c$ and $C_2 = ad - bc.$
Let us choose the symplectic leaf ${\mathcal O}(\lambda)$ corresponding to $C_1 = c =1$ and $C_2 = ad - bc = \lambda:$ it consists of the matrices of the form
\begin{equation}
\label{Leaf}
 L (\zeta, \lambda) =
\Bigl(
\begin{array}{cc}
a &  ad + \lambda - \zeta\\ 1 & d
\end{array}
\Bigl),
\end{equation}
where $\{a, d\} = 1.$

Take $X = {\bf R}^2$ with the standard symplectic structure and consider the embedding
$\phi_{\lambda} : X \rightarrow LGL(2, {\bf R})$ defined by (\ref{Leaf}).
Then the general construction of Theorem 1 gives a multivalued map because
the re-factorisation problem (\ref{factl}) is not uniquely solvable. Indeed
the form (\ref{Leaf}) is invariant under multiplication by the triangular matrices
of the form
\begin{equation}
\label{tri}
 \hat t = \Bigl(
\begin{array}{cc}
1 &  t \\ 0  & 1
\end{array}
\Bigl), 
\end{equation} 
so $$L (\zeta, \lambda)  L (\zeta, \mu) =  \tilde L (\zeta, \lambda) \tilde L (\zeta, \mu),$$
where $\tilde L (\zeta, \lambda) = L (\zeta, \lambda) \hat t, \quad \tilde L (\zeta, \mu) =  \hat t ^{-1} L (\zeta, \mu).$
The product of two symplectic leaves ${\mathcal O}(\lambda)$ and ${\mathcal O}(\mu)$
is actually a three-dimensional Poisson submanifold.

Similarly the product of $N$ symplectic leaves ${\mathcal O}(\lambda_1) {\mathcal O}(\lambda_2) ... {\mathcal O}(\lambda_N)$ is a Poisson submanifold which has dimension $(N+1)$. Let us consider its $N$-dimensional quotient $P(\lambda_1,..., \lambda_N)$ by the action of group $T \approx {\bf R}$ of triangular matrices
of the form (\ref{tri}): $g \rightarrow \hat t  g \hat t ^{-1} .$
Alternatively one can describe $P(\lambda_1,..., \lambda_N)$ as the result of the Hamiltonian reduction of the symplectic manifold $M^{2N} = {\mathcal O}(\lambda_1) \times {\mathcal O}(\lambda_2) \times ...\times {\mathcal O}(\lambda_N)$ by the following Hamiltonian action of the group $T^N$ 
\begin{equation}
\label{act}
L (\zeta, \lambda_i)  \rightarrow \hat t_i L (\zeta, \lambda_i) \hat t_{i+1} ^{-1}.
\end{equation}
The Hamiltonians of the corresponding flows are $H_i = a_{i-1} + d_i.$
Note that although the action is commutative 
the Hamiltonians $H_i$
do not Poisson commute: $\{ H_i, H_{i-1}\} = 1.$

We claim that the reduced space is actually the phase space of the periodic dressing chain.
Indeed the action of $t_i$ has the form
$$d_{i-1} \rightarrow d_{i-1} - t_i, \quad a_{i} \rightarrow a_{i} + t_i$$
while all other variables remain the same. The space of invariant functions
on the quotient space is generated by the functions $G_i = a_i + d_{i-1}, \,\, i=1,...,N.$
The Poisson brackets for these functions $\{G_i, G_{i+1}\} = 1 = -\{G_{i+1}, G_i\}$
are exactly the same as for the variables $g_i$ in the dressing chain.

{\bf Theorem 3.} {\it The phase space of the periodic dressing chain with odd period $N$ is the Hamiltonian reduction of the product of symplectic leaves ${\mathcal O}(\lambda_1) \times {\mathcal O}(\lambda_2) \times ...\times {\mathcal O}(\lambda_N).$ The reduced transfer maps coincide with integrable Adler's dynamical
systems \cite{A}.}

We need only to note that for odd $N$ the quotient space of $M^{2N}$ by the action (\ref{act}) can be identified with the subset ${\mathcal F}(\lambda_1) \times {\mathcal F}(\lambda_2) \times ...\times {\mathcal F}(\lambda_N),$ where ${\mathcal F}(\lambda)$ is the set of matrices of the form
$$ {\mathcal L} =
\Bigl(
\begin{array}{cc}
f &  f^2 + \lambda - \zeta\\ 1 & f
\end{array}
\Bigl),$$
which is nothing else but the Lax matrix for the Adler's map.

\section*{Concluding remarks.}

We have shown that the theory of Lie Poisson group is behind the Poisson properties of the transfer dynamics related to matrix soliton interaction and Adler's maps. There is no doubt that
this is a very general scheme which should be applicable also to many other examples of Yang-Baxter maps (see e.g. \cite{Skl, NY, O, Et, KNY}). For example, the Yang-Baxter map related to geometric crystals \cite{Et, KNY, SV} can be explained using the standard Poisson Lie structure on $GL_N.$
It would be interesting to analyse from this point of view the discrete top \cite{MV}.

\subsection*{ Acknowledgements}

We are grateful to the organisers of the Hayashibara Forum - 2003 in
Oxford-Warwick  (2-6 June 2003) where this work has been initiated.

One of us (APV) would like to acknowledge the support of European research
programmes ENIGMA (contract MRTN-CT-2004-5652) and MISGAM.

\end{document}